\documentclass[a4paper]{article}

\usepackage{amsmath,amssymb,graphicx}
\usepackage{subfigure}

\usepackage{lscape}
\usepackage{mathrsfs}

\usepackage{algorithmic}
\usepackage[ruled]{algorithm}

\usepackage{multirow}
\usepackage{verbatim}
\usepackage{enumerate}
\usepackage{color}

\newtheorem{theorem}{Theorem}
\newtheorem{corollary}{Corollary}

\usepackage{longtable}

\title{ Randomized Iterative Methods with Alternating Projections }

\author{  Hua Xiang$^a$\thanks{Corresponding author (H. Xiang).
E-mail: hxiang@whu.edu.cn. H. Xiang is supported by the National Natural Science Foundation of China (No.11571265, No.11471253) and NSFC-RGC (No.11661161017).
}
\quad Lin Zhang$^b$\thanks{
E-mail: linyz@zju.edu.cn.  L. Zhang is supported by Natural Science Foundation of Zhejiang Province of China (LY17A010027) and National Natural Science Foundation of China (No.11301124). }
\\ \\
\small{$^a$ School of Mathematics and Statistics, Wuhan University, Wuhan, 430072, P.R. China}\\
{\small $^b$ Institute of Mathematics, Hangzhou Dianzi University, Hangzhou 310018, P.R. China.}
}

\begin{document}

\maketitle

\begin{abstract}

We use a unified framework to summarize sixteen randomized iterative methods including Kaczmarz method, coordinate descent method, etc. Some new iterative schemes are given as well. Some relationships with \textsc{mg} and \textsc{ddm} are also discussed. We analyze the convergence properties of the iterative schemes by using the matrix integrals associated with alternating projectors, and demonstrate the convergence behaviors of the randomized iterative methods by numerical examples.

\end{abstract}




\section{Introduction}

Solving a linear system is a basic task in scientific computing.
Given a real matrix $A \in \mathbb{R}^{m\times n}$ and a real vector $b\in \mathbb{R}^m$,
in this paper we consider the following consistent linear system
\begin{equation}\label{eqn:Ax=b}
 Ax = b.
\end{equation}
For a large scale problem, iterative methods are more suitable than  direct methods. Though the classical iterative methods are generally deterministic \cite{Demmel_Book97,GolubLoan_Book2013}, in this paper, we  consider randomized iterative methods, for example, Kaczmarz method, the coordinate descent (\textsc{cd}) method and their block variants. Motivated by work in \cite{GowerRichtarik_SIMAX15}, we give a very general framework to such randomized iterative schemes.

We need two matrices $Y \in \mathbb{R}^{m \times l}$ and $Z \in \mathbb{R}^{n \times l}$, where $l \leqslant \min \{m,n\}$, and usually $l$ is much smaller than $m$ and $n$. In most practical cases, only one of  the matrices $Y$ and $Z$ is given, and the other one is determined consequently. Suppose that $Y$ and $Z$ are of full column rank. Let $ E \equiv Y^T A Z$, which is invertible in most cases, otherwise we use the pseudoinverse in the following. For convenience, we define
$$ \Xi := Z E^{\dag} Y^T = Z ( Y^T A Z )^{\dag} Y^T, $$
where $\dag$ denotes the (Moore-Penrose) pseudoinverse.

In this paper, we consider the following iterative scheme:
\begin{equation} \label{eqn:GeneralIterationXi}
 x^{k+1} = x^{k} + \Xi ( b - A x^k ),
\end{equation}
where $x^k$ ($k=0,1,\cdots$) is the $k$th approximate solution, and $x^0$ is the initial guess.
This is our general framework under consideration.
Let $ T := I-\Xi A= I - Z ( Y^T A Z )^{\dag}Y^T A $. Then \eqref{eqn:GeneralIterationXi} can be rewritten as
\begin{equation} \label{eqn:GeneralIterationT}
 x^{k+1} = T x^{k} + \Xi b .
\end{equation}
Assume that $x^*$ is the true solution of \eqref{eqn:Ax=b}, and define the $k$th iteration error $e^k = x^k - x^*$. Then   we have the error propagation:
\begin{equation}
 e^{k+1} = ( I - \Xi A ) e^k = T e^k.
\end{equation}

Besides the formulations \eqref{eqn:GeneralIterationXi} and \eqref{eqn:GeneralIterationT}, the iteration scheme can be viewed in seemingly different but equivalent ways as pointed out in \cite{GowerRichtarik_SIMAX15}. From an algebraic viewpoint, $x^{k+1}$ is the solution of $Y^T A x =  Y^T b$ and $x = x^k + Zy$. From a geometric viewpoint, $x^{k+1}$ is the random intersect of two affine spaces. That is, $ \{ x^{k+1} \} = ( x^* + \text{Null}(Y^T A) ) ~ \cap ~ (x^k + \text{Range}(Z))$.

We can easily check that $\Xi A$ and $T$ are projectors. If $Y$ or $Z$ is different in each step $k$, then the scheme is nonstationary and $T$ depends on $k$, denoted by $T_k$.
 Choosing $Y$ or $Z$ properly in each step, we can recover the classical Kaczmarz method, and the coordinate descent method, etc. We will show how to achieve this in Section 2.
Furthermore, from the error propagation, we have
\begin{equation}\label{eqn:errorPropagation}
e^{k+1} = T_k \cdots T_0 e^0,
\end{equation}
where $T_0, \cdots, T_k$ are a series of different alternating projectors applied to the initial iteration error $e^0$. We show that under such alternating projections the final iteration error approaches zero in the sense of probability in Section 3. The convergence is proved by using the matrix integrals associated projections, without the assumption on the matrix positive definiteness in \cite{GowerRichtarik_SIMAX15}.  Numerical examples are given in Section 4.

\section{Iterative schemes}

We consider the iterative methods under the unified framework of \eqref{eqn:GeneralIterationXi} or \eqref{eqn:GeneralIterationT}. Depending on the choice of parameter matrices $Y$ and $Z$, we divide the iterative schemes into three types for convenience. In the following we always assume that the parameter matrix $G$ is symmetric positive definite (\textsc{spd}).

(1) Type-\textsc{k} methods. In this type, $Y$ is given. We set $Z= GA^T Y$, then consequently 
$\Xi A = Z E^{\dag} Y^T A = GA^T Y(Y^TAGA^TY)^{\dag} Y^TA  $. These are row action methods, and the Kaczmarz algorithm is the typical method of this type.

For this type, as pointed in \cite{GowerRichtarik_SIMAX15}, the iterative scheme can be expressed from a sketching viewpoint as
 $$x^{k+1} = \arg \min ||x-x^k||_{G^{-1}}, ~ \text{s.t.},~ Y^T Ax= Y^T b. $$
And also it can be rewritten from the optimization viewpoint as follows.
 $$x^{k+1} = \arg \min ||x-x^*||_{G^{-1}}, ~ \text{s.t.},~   x = x^k + Z y. $$

(2) Type-\textsc{c} methods. These are column action methods, and the coordinate descent algorithm is the typical method of this type. In this type, $Z$ is given. we set $Y=GA Z$, then 
$\Xi A = Z E^{\dag} Y^T A = Z(Z^TA^T G AZ)^{\dag} Z^TA^T G  A $.

Similarly, for this type methods we can reexpress the iterative scheme from a sketching viewpoint as
$$x^{k+1} = \arg \min || x-x^k||_{\widehat{G}}, ~ \text{s.t.},~  Y^T  Ax= Y^T  b, $$
where $A$ is supposed to be of full column rank, and then $\widehat{G} := A^T G A$ is \textsc{spd}.
And also we can reformulate from the optimization viewpoint
$$ x^{k+1}= \arg \min || x-x^* ||_{\widehat{G}}, ~ \text{s.t.}, ~ x=x^k + Z y. $$

(3) Type-\textsc{s} methods. This type of iterative schemes is applied to the symmetric cases. In this type, we choose $Y=Z$, then $E=Z^T A Z$, $\Xi=Z(Z^TAZ)^{\dag} Z^T$, and $\Xi A = Z(Z^TAZ)^{\dag} Z^T A$. For the \textsc{spd} case, under the settting $G=A^{-1}$, type-\textsc{k} and type-\textsc{c} methods will be the same, and both reduce to type-\textsc{s} methods.

In the following, by choosing specific $Y$ or $Z$, we recover the classical well-known iteration schemes, and also obtain some new ones. The iteration schemes are summarize in Tables \ref{tab:KCDseriesMethods} and \ref{tab:KCDseriesSchemes}. We will derive and discuss them one by one.

\begin{table}[!htbp]
\begin{center}
\caption{\label{tab:KCDseriesMethods} Matrices in the randomized iterative methods}
\begin{tabular}{llllll}\hline
No. & $Y$ &  $Z$ & $E=Y^TAZ$ & $\Xi=Z E^{\dag} Y^T$ & $T=I-\Xi A$  \\ \hline
K1 & $e_i$ & $A^T e_i$ & $e_i^T A A^T e_i$ & $ A^T e_i (e_i^T A A^T e_i)^{\dag} e_i^T$ & $I - A^T e_i (e_i^T A A^T e_i)^{\dag} (A^T e_i)^T$  \\
K2 & $\omega$ & $ A^T \omega$ & $\omega^T A A^T \omega $ & $A^T \omega (\omega^T A A^T \omega)^{\dag} \omega^T $ & $ I - A^T \omega (\omega^T A A^T \omega)^{\dag} \omega^T A$ \\
K3 & $I_R $ & $ A^T I_R$ & $I_R^T A A^T I_R $ & $ A^T I_R (I_R^T A A^T I_R)^{\dag} I_R^T$ & $ I - A^T I_R (I_R^T A A^T I_R)^{\dag} I_R^T A$  \\
K4 & $\Omega $ & $ A^T \Omega$ & $ \Omega^T A A^T \Omega$ & $ A^T \Omega (\Omega^T A A^T \Omega)^{\dag} \Omega^T$ & $I - A^T \Omega (\Omega^T A A^T \Omega)^{\dag} \Omega^T A $   \\
K5 & $ I_R$ & $ G A^T I_R$ & $I_R^T A G A^T I_R $ & $ G A^T I_R (I_R^T A G A^T I_R)^{\dag} I_R^T$ & $I - G A^T I_R (I_R^T A G A^T I_R)^{\dag} I_R^T A $   \\
K6 & $ \Omega$ & $ G A^T \Omega$ & $ \Omega^T A G A^T \Omega$ & $ G A^T \Omega (\Omega^T A G A^T \Omega)^{\dag} \Omega^T$ & $ I - G A^T \Omega (\Omega^T A G A^T \Omega)^{\dag} \Omega^T  A $  \\
C1 & $ A e_j$ & $e_j $ & $ e_j^T A^T  A e_j$ & $e_j (e_j^T A^T A e_j)^{\dag} (A e_j)^T $ & $ I - e_j (e_j^T A^T A e_j)^{\dag} (A e_j)^T A$ \\
C2 & $A \omega $ & $ \omega $ & $ \omega^T A^T A \omega$ & $ \omega (\omega^T A^T A \omega)^{\dag} \omega^T A^T$ & $ I - \omega (\omega^T A^T A \omega)^{\dag} \omega^T A^T  A $ \\
C3 & $ A I_C$ & $ I_C$ & $I_C^T A^T  A I_C $ & $ I_C ( I_C^T A^T A I_C )^{\dag} I_C^T A^T$ & $ I - I_C ( I_C^T A^T A I_C )^{\dag} I_C^T A^T A$ \\
C4 & $A \Omega $ & $ \Omega$ & $\Omega^T A^T A \Omega $ & $ \Omega (\Omega^T A^T A \Omega)^{\dag} \Omega^T A^T$ & $I - \Omega (\Omega^T A^T A \Omega)^{\dag} \Omega^T A^T A $ \\
C5 & $G A I_C $ & $I_C $ & $ I_C^T A^T G A I_C$ & $I_C ( I_C^T A^T G A I_C )^{\dag} I_C^T A^T G $ & $ I - I_C ( I_C^T A^T G A I_C )^{\dag} I_C^T A^T G A $  \\
C6 & $G A \Omega $ & $\Omega $ & $\Omega^T A^T G A \Omega $ & $ \Omega ( \Omega^T A^T G A \Omega )^{\dag} \Omega^T A^T G$ & $ I - \Omega ( \Omega^T A^T G A \Omega )^{\dag} \Omega^T A^T G A $ \\
S1 & $ e_i$ & $ e_i $ & $ e_i^T  A e_i$ & $ e_i (e_i^T A e_i)^{\dag} e_i^T$ & $I - e_i (e_i^T A e_i)^{\dag} e_i^T A  $  \\
S2 & $ \omega$ & $ \omega$ & $ \omega^T  A \omega$ & $ \omega (\omega^T  A \omega)^{\dag} \omega^T$ & $ I - \omega (\omega^T  A \omega)^{\dag} \omega^T  A$  \\
S3 & $ I_C$ & $ I_C$ & $ I_C^T   A I_C$ & $ I_C ( I_C^T  A I_C )^{\dag} I_C^T $ & $ I - I_C ( I_C^T  A I_C )^{\dag} I_C^T  A$ \\
S4 & $ \Omega$ & $ \Omega$ & $ \Omega^T  A \Omega $ & $ \Omega (\Omega^T  A \Omega)^{\dag} \Omega^T$ & $ I - \Omega (\Omega^T  A \Omega)^{\dag} \Omega^T A$  \\
\hline
\end{tabular}
\end{center}
\end{table}


\begin{table}[!htbp]
\begin{center}
\caption{\label{tab:KCDseriesSchemes} Randomized iterative schemes}
\begin{tabular}{lll}\hline
No. & Schemes & Remarks\\ \hline
K1 & $ x^{k+1} = x^k + \frac{b_i - A_{i:} x^k}{||A_{i:}||^2}(A_{i:})^T. $ & Randomized Kaczmarz\\
K2 & $x^{k+1} = x^k + \frac{\omega^T (b - A x^k) }{||A^T \omega||^2} A^T \omega$ & Gaussian Kaczmarz\\
\multirow{2}{*}{K3} & \multirow{2}{*}{$x^{k+1} = x^k + A_{R:}^T (A_{R:} A_{R:}^T)^{\dag} ( b_R - A_{R:} x^k )$} & Randomized block Kaczmarz; \\
   &  & \textsc{psh} methods, see (83) in \cite{Maess_JCAM88} \\
K4 & $x^{k+1} = x^k + A^T \Omega (\Omega^T A A^T \Omega)^{\dag} \Omega^T (b-Ax^k)$ & see (88) in \cite{Maess_JCAM88}, but not stochastic \\
K5 & $x^{k+1} = x^k +  G A_{R:}^T (A_{R:} G A_{R:}^T)^{\dag}  (b_R- A_{R:} x^k)$ & see (83') in \cite{Maess_JCAM88}  \\
K6 & $x^{k+1} = x^k + G A^T \Omega (\Omega^T A G A^T \Omega)^{\dag} \Omega^T (b-Ax^k)$ & see (2.8) in \cite{GowerRichtarik_SIMAX15}  \\
\multirow{2}{*}{C1} & \multirow{2}{*}{$x^{k+1} = x^k + \frac{A_{:j}^T(b-Ax^k )}{||A_{:j}||^2} e_j$} & Randomized coordinate descent, \textsc{rgs} \\ 
   &  & De la Garza column algorithm (1951) \cite{Maess_JCAM88} \\
C2 & $x^{k+1} = x^k + \frac{\omega^T A^T (b - A x^k) }{||A \omega||^2}  \omega $ & Gaussian \textsc{ls} \cite{GowerRichtarik_SIMAX15}\\
C3 & $x^{k+1} = x^k + I_C ( A_{:C}^T A_{:C} )^{\dag} A_{:C}^T (b-A x^k) $  & \textsc{spa} methods, see (87) in \cite{Maess_JCAM88}  \\
C4 & $x^{k+1} = x^k +  \Omega (\Omega^T A^T A \Omega)^{\dag} \Omega^T A^T (b - A x^k)$ & It's new to the best of our knowledge \\
C5 & $x^{k+1} = x^k +  I_C ( A_{:C}^T G A_{:C} )^{\dag} A_{:C}^T G (b - A x^k)$ & see (87') in \cite{Maess_JCAM88} \\
C6 & $x^{k+1} = x^k +  \Omega ( \Omega^T A^T G A \Omega )^{\dag} \Omega^T A^T G (b - A x^k)$ & It's new to the best of our knowledge \\
S1 & $x^{k+1} = x^k + \frac{ b_i - A_{i:}^T x^k }{A_{ii}} e_i $  & Randomized \textsc{cd}: positive definite \\
S2 & $x^{k+1} = x^k + \frac{\omega^T (b - A x^k) }{|| \omega||_A^2}  \omega $ & Gaussian \textsc{cd}: positive definite\\
S3 & $x^{k+1} = x^k + I_{:C} ( I_{:C}^T A I_{:C} )^{\dag} I_{:C}^T (b-A x^k)$  & Randomized Newton: positive definite \\
S4 & $x^{k+1} = x^k +  \Omega (\Omega^T  A \Omega)^{\dag} \Omega^T (b - A x^k)$  & It's new to the best of our knowledge \\
\hline
\end{tabular}
\end{center}
\end{table}

\subsection{Kaczmarz and row methods}


The Kaczmarz method, also known as the algebraic reconstruction technique (\textsc{art}), is a typical row action method, dated back to the Polish mathematician Stefan Kaczmarz in 1937 and reignited by Strohmer and Vershynin \cite{StrohmerVershynin_JFAA09}.  
The advantage of Kaczmarz method and the following \textsc{cd} lies in the fact that at a time they only need access to individual rows (or columns) rather than the entire coefficient matrix.
Due to its simplicity, it has numerous applications in image reconstruction, signal processing, etc.

In the Kaczmarz algorithm, the approximate solution $x$ is projected onto the hyperplane determined by a random row of the matrix $A$ and the respective element of the right hand side $b$.
Let $A^T = (\alpha_1, \cdots, \alpha_m)$. That is, $\alpha_i^T$ is the $i$th row of $A$.  The $i$th component of \eqref{eqn:Ax=b} reads $\alpha_i^T x = b_i$, i.e., $\frac{\alpha_i^T x}{|| \alpha_i||} = \frac{ b_i}{||\alpha_i||} (i=1,\cdots, m) $.
For an approximation $x$ not in the plane, the distance of this point $x$ to the plane is $ \left( \frac{ b_i}{||\alpha_i||} - \frac{\alpha_i^T x}{|| \alpha_i||} \right) \frac{\alpha_i}{|| \alpha_i||} $, which is used to update the approximate solution. Then the iterative scheme reads:
\begin{equation}\label{eqn:Kaczmarz1}
 x \leftarrow x +   \frac{ b_i-\alpha_i^T x}{|| \alpha_i||^2}  \alpha_i ,
\end{equation}
which is in fact achieved by solving $\arg\min_{ \{z | \alpha_i^T z = b_i\} } \frac12 ||z- x||^2$.
This iterative scheme is also equivalent to Gauss-Seidel on $ AA^Ty=b$ with the standard primal-dual mapping $x=A^Ty$.

We will show that Kaczmarz method is a special case of the general iteration scheme \eqref{eqn:GeneralIterationXi}. Let $ Y = e_i$, where $e_i$ is the $i$th column of the unit matrix $I_m$, and set $Z=A^T e_i$, then $E = e_i^T A A^T e_i$, $\Xi = Z E^{\dag} Y^T = A^T e_i (e_i^T A A^T e_i)^{\dag} e_i^T$.
Then the iteration scheme \eqref{eqn:GeneralIterationXi} becomes
\begin{equation}\label{eqn:tagK1}
 x^{k+1} = x^k + \frac{b_i - A_{i:} x^k}{||A_{i:}||^2}(A_{i:})^T,  \tag{K1}
\end{equation}
which is in fact the same as \eqref{eqn:Kaczmarz1}, where $A_{i:}$ denotes the $i$th row of $A$, similar to the \textsc{Matlab} notation.
The related matrices and the iterative schemes are given in Tables \ref{tab:KCDseriesMethods} and \ref{tab:KCDseriesSchemes} respectively, where the Kaczmarz method is denoted by K1. Note that the corresponding iteration matrix $T = I - (A^T e_i) (A^T e_i)^T / ||A^T e_i||_2^2 $ is an alternating projector, depending on $i$ and changing in each iteration step.

Relaxation parameters can be introduced in \eqref{eqn:tagK1}, and  it can be also extended to nonlinear versions.
Combining with the soft shrinkage, the Kaczmarz method can be used for sparse solutions \cite{LorenzSchopferWenger_SIIS14}.
It aslo appears as a special case of the current popular stochastic gradient descent (\textsc{sgd}) method for convex optimization.


Instead of the choice $Y=e_i$ which corresponds to  randomly choosing one row in each iteration of K1, we can extend to the real Gaussian vector, that is, $Y = \omega$, where $\omega = (\omega_1, \cdots, \omega_m)^T$, and $\omega_j \sim N(0,1)$, the standard normal distribution. Let $Z = A^T \omega$, then $ \Xi = A^T \omega (\omega^T A A^T \omega)^{\dag} \omega^T$. We obtain the Gaussian Kaczmarz method (K2 in Tables \ref{tab:KCDseriesMethods} and \ref{tab:KCDseriesSchemes}):
\begin{equation}\label{eqn:tagK2}
x^{k+1} = x^k + \frac{\omega^T (b - A x^k) }{||A^T \omega||^2} A^T \omega .  \tag{K2}
\end{equation}


The above two methods can be extended to the block variants. Instead of just 
 using one column vector $e_i$ or $\omega$, 
 we can work on several columns simultaneously. Let $R$ be a  random subset including row indices, and correspondingly $I_R$, a column concatenation of the columns of $I$ indexed by $R$. We extend the K1 method by defining $ Y = I_R$. And let $ Z=A^T I_R$, then $\Xi = A^T I_R (I_R^T A A^T I_R)^{\dag} I_R^T$, and the randomized block Kaczmarz iteration scheme (K3 in Tables \ref{tab:KCDseriesMethods} and \ref{tab:KCDseriesSchemes}) reads
\begin{equation}\label{eqn:tagK3}
 x^{k+1} = x^k + A_{R:}^T (A_{R:} A_{R:}^T)^{\dag} ( b_R - A_{R:} x^k ), \tag{K3}
\end{equation}
where  $A_{R:}$ is formed by the rows of $A$ indexed by $R$, and $A_{R:}^T $ stands for $( A_{R:} )^T $.


The block extension of \eqref{eqn:tagK2} is natural, just letting $Y$ consist of  several columns of Gaussian vectors. That is, $ Y = \Omega$, where $\Omega$ is a Gaussian matrix with i.i.d. entries. Obviously, $\Omega$ and $I_R$ are the generalizaiton of $\omega$ and $e_i$ in \ref{eqn:tagK2} and \ref{eqn:tagK1}, respectively.
Correspondingly define $Z = A^T \Omega$, then we have the iteration scheme (K4 in Tables \ref{tab:KCDseriesMethods} and \ref{tab:KCDseriesSchemes})
\begin{equation}\label{eqn:tagK4}
 x^{k+1} =    x^k + A^T \Omega (\Omega^T A A^T \Omega)^{\dag} \Omega^T (b-Ax^k). \tag{K4}
\end{equation}


If we introduce a symmetric positive definite (\textsc{spd}) matrix $G$, and use an energetic norm $||x||_G := \sqrt{x^T G x}$ instead of the Euclidean norm in \eqref{eqn:tagK1}, we can obtain more general iterative scheme (see \cite[formula (83')]{Maess_JCAM88} and the references therein). In this way we can generalize  the methods \ref{eqn:tagK3} and \ref{eqn:tagK4}, where the matrix $Y$ is unchanged, but the matrix $Z$ is modified by $G$ (see below). We describe how to achieve such extensions in the following.

Choose $Y = I_R$ and define $Z= G A^T I_R$, then $\Xi = G A^T I_R (I_R^T A G A^T I_R)^{\dag} I_R^T$, and the iteration scheme (K5 in Tables \ref{tab:KCDseriesMethods} and \ref{tab:KCDseriesSchemes}) reads
\begin{equation}\label{eqn:tagK5}
x^{k+1} = x^k +  G A_{R:}^T (A_{R:} G A_{R:}^T)^{\dag}  (b_R- A_{R:} x^k). \tag{K5}
\end{equation}


Choose $Y = \Omega$ and define $Z= G A^T \Omega$, then $\Xi = G A^T \Omega (\Omega^T A G A^T \Omega)^{\dag} \Omega^T$, and the iteration scheme (K6 in Tables \ref{tab:KCDseriesMethods} and \ref{tab:KCDseriesSchemes}) reads
\begin{equation}\label{eqn:tagK6}
 x^{k+1}  = x^k + G A^T \Omega (\Omega^T A G A^T \Omega)^{\dag} \Omega^T (b-Ax^k), \tag{K6}
\end{equation}
which has been intensively investigated in \cite{GowerRichtarik_SIMAX15}.

\subsection{Coordinate descent and column methods}

Randomized coordinate descent (\textsc{cd}) algorithm, also known as the randomized Gauss-Seidel (\textsc{rgs}), is another simple and popular method for linear systems and has been around for a long time \cite{Zangwill_Book69}.   It use gradient information $\nabla \phi$, where $\phi(x) = ||Ax-b||^2$, about a single coordinate to update the approximate solution in each iteration. Besides, it can reduce to the Kaczmarz method when working on $AA^T y=b$ with $x=A^Ty$.
Under our framework \eqref{eqn:GeneralIterationXi},  the difference between Kaczmarz and \textsc{cd} arises from the different choices on the matrices $Y$ and $Z$. For such column methods in this subsection, we first choose $Z$, and then set $Y=AZ$.

Choose $Z= e_j$, and set $Y = A e_j$, then we have $ \Xi = Z E^{\dag} Y^T  = e_j (e_j^T A^T A e_j)^{\dag} (A e_j)^T $, and the iteration scheme (C1 in Tables \ref{tab:KCDseriesMethods} and \ref{tab:KCDseriesSchemes})
\begin{equation}\label{eqn:tagC1}
 x^{k+1} = x^k + \frac{A_{:j}^T(b-Ax^k )}{||A_{:j}||^2} e_j , \tag{C1}
\end{equation}
which is the randomized coordinate descent method, with $j$ being a random number from $\{1,\cdots, n\}$.
At iteration $j$, the current residual is projected randomly onto a column of the matrix $A$.
Such randomized updating needs small costs per iteration and gives provable global convergence. Before this, various strategies have been considered for picking the coordinate to update, such as cyclic coordinate update and the best coordinate update, however these schemes are either hard to estimate or difficult to be implemented efficiently.


Instead of just choosing one column, we choose $Z = \omega \sim N(0,1)$, and set $ Y = A \omega$, a linear combination of the columns of $A$, then $\Xi =  \omega (\omega^T A^T A \omega)^{\dag} \omega^T A^T $, and the Gaussian \textsc{ls}  iteration scheme \cite{GowerRichtarik_SIMAX15} (C2 in Tables \ref{tab:KCDseriesMethods} and \ref{tab:KCDseriesSchemes}) reads
\begin{equation}\label{eqn:tagC2}
x^{k+1} = x^k + \frac{\omega^T A^T (b - A x^k) }{||A \omega||^2}  \omega . \tag{C2}
\end{equation}


We extend the methods \ref{eqn:tagC1} and \ref{eqn:tagC2} to their block variants.
Let $C$ be a random subset containing column indices, and correspondingly $I_C$ be the column concatenation of the columns of $I_n$ indexed by $C$. Choose $Z= I_C$, and set $Y = A I_C$. We can check that $\Xi = I_C ( I_C^T A^T A I_C )^{\dag} I_C^T A^T$, and have the following scheme (C3 in Tables \ref{tab:KCDseriesMethods} and \ref{tab:KCDseriesSchemes})
\begin{equation}\label{eqn:tagC3}
 x^{k+1} = x^k + I_C ( A_{:C}^T A_{:C} )^{\dag} A_{:C}^T (b-A x^k)  . \tag{C3}
\end{equation}
Obviously this is the block version of \ref{eqn:tagC1}.


Replacing the Gaussian vector $\omega$ in \ref{eqn:tagC2} by the Gaussian matrix $\Omega$ with i.i.d Gaussian normal entries, that is, setting $Z = \Omega$ and $Y = A \Omega$, we have the block version of \ref{eqn:tagC2}. The resulting scheme (C4 in Tables \ref{tab:KCDseriesMethods} and \ref{tab:KCDseriesSchemes}) reads
\begin{equation}\label{eqn:tagC4}
 x^{k+1} =  x^k + \Xi ( b-Ax^k) = x^k +  \Omega (\Omega^T A^T A \Omega)^{\dag} \Omega^T A^T (b - A x^k). \tag{C4}
\end{equation}
To the best of our knowledge, the scheme \eqref{eqn:tagC4} is new.


Similar to the schemes \ref{eqn:tagK5} and \ref{eqn:tagK6}, we extend the methods \ref{eqn:tagC3} and \ref{eqn:tagC4} by using the \textsc{spd} matrix $G$, where the matrix $Z$ is unchanged, but the matrix $Y$ is modified by $G$, compared with those matrices in \ref{eqn:tagC3} and \ref{eqn:tagC4}. The extensions are described as follows.

Choose $Z= I_C$ and define $Y = G A I_C$, then $ \Xi = I_C ( I_C^T A^T G A I_C )^{\dag} I_C^T A^T G$,  and the iteration scheme (C5 in Tables \ref{tab:KCDseriesMethods} and \ref{tab:KCDseriesSchemes}) reads
\begin{equation}\label{eqn:tagC5}
 x^{k+1} =   x^k +  I_C ( A_{:C}^T G A_{:C} )^{\dag} A_{:C}^T G (b - A x^k). \tag{C5}
\end{equation}


Choose $Z= \Omega$ and define $Y = G A \Omega$, then $\Xi = \Omega ( \Omega^T A^T G A \Omega )^{\dag} \Omega^T A^T G$, and the iteration scheme (C6 in Tables \ref{tab:KCDseriesMethods} and \ref{tab:KCDseriesSchemes}) reads
\begin{equation}\label{eqn:tagC6}
 x^{k+1} =  x^k + \Xi ( b-Ax^k) = x^k +  \Omega ( \Omega^T A^T G A \Omega )^{\dag} \Omega^T A^T G (b - A x^k), \tag{C6}
\end{equation}
which is new to the best of our knowledge.

\subsection{Schemes for symmetric cases}


For the \textsc{spd} matrix $A$, it is reasonable to choose $Y=Z$ to keep the matrix $E=Y^T A Z$ is \textsc{spd} consequently. In the following, we choose specific matrices $Y$ and derive four different iterative schemes.

Similar to \ref{eqn:tagK1} and \ref{eqn:tagC1}, we choose $Y = Z = e_i$. Then we have $ \Xi = e_i (e_i^T A e_i)^{\dag} e_i^T $, and the iteration scheme (S1 in Tables \ref{tab:KCDseriesMethods} and \ref{tab:KCDseriesSchemes})
\begin{equation}\label{eqn:tagS1}
x^{k+1} = x^k + \frac{ b_i - A_{i:}^T x^k }{A_{ii}} e_i ,  \tag{S1}
\end{equation}
where $A_{i:} = A_{:i}^T$.


Extending the schemes \ref{eqn:tagK2} and \ref{eqn:tagC2} to the \textsc{spd} case, we choose $ Y = Z = \omega$. We can verify that $\Xi =  \omega (\omega^T  A \omega)^{\dag} \omega^T $, and the corresponding iteration scheme (S2 in Tables \ref{tab:KCDseriesMethods} and \ref{tab:KCDseriesSchemes}) reads
\begin{equation}\label{eqn:tagS2}
 x^{k+1} = x^k + \frac{\omega^T (b - A x^k) }{|| \omega||_A^2}  \omega . \tag{S2}
\end{equation}


The schemes \ref{eqn:tagS1} and \ref{eqn:tagS2} can be generalized to the block versions. Extending the schemes \ref{eqn:tagK3} and \ref{eqn:tagC3} to the \textsc{spd} case by choosing $ Y = Z = I_C$, we have $\Xi = I_C ( I_C^T  A I_C )^{\dag} I_C^T $, and the iteration scheme (S3 in Tables \ref{tab:KCDseriesMethods} and \ref{tab:KCDseriesSchemes})
\begin{equation}\label{eqn:tagS3}
x^{k+1} = x^k + I_{:C} ( I_{:C}^T A I_{:C} )^{\dag} I_{:C}^T (b-A x^k)  . \tag{S3}
\end{equation}


Extending the schemes \ref{eqn:tagK4} and \ref{eqn:tagC4} to the \textsc{spd} case by choosing $ Y = Z = \Omega$, we obtain $ \Xi = \Omega (\Omega^T  A \Omega)^{\dag} \Omega^T $, and the iteration scheme (S4 in Tables \ref{tab:KCDseriesMethods} and \ref{tab:KCDseriesSchemes})
\begin{equation}\label{eqn:tagS4}
x^{k+1} =  x^k + \Xi ( b-Ax^k) = x^k +  \Omega (\Omega^T  A \Omega)^{\dag} \Omega^T (b - A x^k). \tag{S4}
\end{equation}
To the knowledge of the authors, the scheme \ref{eqn:tagS4} is new.

We can similarly introduce the matrix $G$, and choose $Y =Z= G I_C$ or $Y =Z = G \Omega$, but the resulting iterative schemes will not be as natural as the extensions such as K5-6 and C5-6, and hence we omit such discussion here. Besides, for the \textsc{spd} case, if we choose $G=A^{-1}$ in K5-6 and C5-6, then K5 and C5 are the same and reduce to S3, while K6 and C6 are the same and reduce to S4.

\subsection{Related matrices in MG \& DDM}

We can find the related matrices in the multigrid (\textsc{mg}) method and the domain decomposition method (\textsc{ddm}).
In \textsc{mg} method, $Z$ and $Y$ denote a prolongation and a restriction operator respectively. In \textsc{ddm}, $Z$ and $Y$ are full rank matrices which span the coarse grid subspaces. $E= Y^T A Z$ is the Galerkin matrix or coarse-grid matrix, and $ \Xi = Z E^{-1} Y^T $ is the coarse-grid correction matrix (see \cite{NatafXiangDoleanSpillane_SISC11} and the references therein).

In \textsc{ddm}, the balancing Neumann-Neumann preconditioner  and the \textsc{feti} algorithm have been intensively investigated, see
\cite{ToselliWidlund_book2005} and references therein. For symmetric
systems the balancing preconditioner was proposed by Mandel
\cite{Mandel_CANM93}. For nonsymmetric systems the abstract balancing preconditioner reads \cite{ErlanggaNabben_SIMAX08},
\begin{equation*}
 P_\textsc{bnn} := T B^{-1} S + Z E^{-1} Y^T,
\end{equation*}
where $T = I -  \Xi A$, $S = I - A \Xi$, and $B$ is a preconditioner for $A$.
We can check that the
error transfer operator corresponding to $P_\textsc{bnn}$ is
\begin{equation}
T_\textsc{bnn}:= I-P_\textsc{bnn} A = T (I-B^{-1}A) T,
\end{equation}
where $T = S^T$ for the symmetric case.

The multigrid V(1,1)-cycle preconditioner $P_\textsc{mg}$ with the smoother $M$ is explicitly given by \cite{TangMaclachlanNabbenVuik_SIMAX10}
$$ P_\textsc{mg}:= M^{-T} T + T^T M^{-1} +  Z E^{-1} Z^T - M^{-T} T A M^{-1}.$$
The error propagation matrix of \textsc{mg} V(1,1)-cycle
preconditioner reads
\begin{equation}
T_\textsc{mg}:=I-P_\textsc{mg}A = (I-
M^{-T} A ) T (I - M^{-1} A).
\end{equation}

Choosing the proper smoother $B$ in $P_\textsc{bnn}$
one can ensure that
$P_\textsc{mg}$ and $P_\textsc{bnn}$ are \textsc{spd}, and  that $P_\textsc{mg} A$ and $P_\textsc{bnn} A$   have the
same spectrum \cite{TangMaclachlanNabbenVuik_SIMAX10}.
The difference lies in the fact that the smoother is used two
times in \textsc{mg} while the coarse grid correction is applied two times in \textsc{ddm}.

Besides, the deflation technique is another closely related method.
The deflation can be used in the following way to solve the linear
system \eqref{eqn:Ax=b}. The solution $x$ of \eqref{eqn:Ax=b} is decomposed into two parts \cite{ErlanggaNabben_SIMAX08}
$$x =
(I-T) x + T x=  \Xi  b + T x. $$
 The first term in the last formula can be easily
computed, and the second term is obtained via solving a singular linear
system,
\begin{equation*}
S A x  = S b,
\end{equation*}
which is consistent and solvable by  applying a Krylov subspace method
for nonsymmetric systems, e.g., \textsc{gmres} or \textsc{b}i\textsc{cg}stab. Its solution   is
non-unique, but $T x$ is unique. 
The corresponding deflated preconditioning system reads $ B^{-1} S
A y = B^{-1} S b  $,
where $B^{-1} S A$ and $P_\textsc{bnn} A$  have the same spectra  except that the
zero eigenvalues of $B^{-1} S A$ are shifted to ones in $P_\textsc{bnn} A$
\cite{ErlanggaNabben_SIMAX08}.

\section{Convergence analysis}

Firstly, we examine the convergence by using the spectra information, especially the minimum eigenvalue of the associated matrix. The convergence of type-\textsc{k} methods was already proved in \cite{GowerRichtarik_SIMAX15}. We present the convergence of type-\textsc{c} and type-\textsc{s} methods, but the techniques are the same as those in \cite{GowerRichtarik_SIMAX15}. Secondly, we give a new and unified proof these iterative methods by using the average properties of random alternating projectors through matrix integral \cite{Zhang_JPAMT17}, without the assumption about positive definiteness used in \cite{GowerRichtarik_SIMAX15}, where such assumption relates to the properties of the random sampling matrices and the coefficient matrix. But our convergence proof still needs the assumption that the coefficient matrix is of full column rank.

\subsection{Convergence results using spectra information}
The convergence of Kaczmarz method \eqref{eqn:tagK1} was considered in \cite{StrohmerVershynin_JFAA09}.
Suppose that the index $i$ is chosen with probability proportional to the magnitude of the $i$th row of $A$. It results in a  convergence with \cite{StrohmerVershynin_JFAA09}
$$\mathbb{E}[||e^k||_2^2] \leqslant \left( 1-\frac{\lambda_{\min}(A^TA)}{||A||_F^2} \right)^k ||e^0||_2^2. $$
The randomized Kaczmarz converges only when the system \eqref{eqn:Ax=b} is consistent. Otherwise, for noisy linear systems it hovers around the least squares (\textsc{ls}) solution within guaranteed bounds \cite{Needell_BIT10}.
The extended Kaczmarz (\textsc{ek}) algorithm \cite{ZouziasFreris_SIMAX13} consisting of \textsc{cd} and Kaczmarz iterations can fix this drawback and converge to the \textsc{ls} solution.

The convergence of coordinate descent method \eqref{eqn:tagC1} was considered in \cite{LeventhalLewis_MOR10}.
Suppose that the index $j$ is chosen with probability proportional to the magnitude of the $i$th column of $A$. It results in a  convergence with \cite{LeventhalLewis_MOR10}
$$\mathbb{E}[||e^k||_{A^TA}^2] \leqslant \left( 1-\frac{\lambda_{\min}(A^TA)}{||A||_F^2} \right)^k ||e^0||_{A^TA}^2. $$
The column methods like \textsc{cd} compute a \textsc{ls} solution, unlike the row methods which exhibit cyclic convergence and aim to a minimum-norm solution to a consistent systems.

The convergence of the scheme \eqref{eqn:tagS1} was also investigated in \cite{LeventhalLewis_MOR10}, where $A$ is \textsc{spd}, and $i$ in \eqref{eqn:tagS1} is chosen randomly according to the probability distribution $p_i = A_{ii} / \mathrm{trace}(A)$. Then the convergence result reads \cite{LeventhalLewis_MOR10}
$$ \mathbb{E} [ ||e^k||_A^2 ] \leqslant \left( 1- \frac{\lambda_{\min}(A)}{\mathrm{trace}(A)} \right)^k ||e^0||_A^2. $$
More convergence results can be found in \cite{GowerRichtarik_SIMAX15} and the references therein.
Define $V:= A^T \Omega (\Omega^T A G A^T \Omega )^{\dag} \Omega^T A $, then $\Xi A = G V$ and $T=I-\Xi A =I-GV$, where $G$ is \textsc{spd}.
A general convergence for \ref{eqn:tagK6} was given in \cite{GowerRichtarik_SIMAX15} under the assumption that $\mathbb{E}[V]$ is positive definite. The   convergence result of the scheme \eqref{eqn:tagK6} reads
\begin{equation}\label{eqn:K6convergence}
 \mathbb{E}[ ||e^k||_{G^{-1}}^2 \leqslant  \rho^k   ||e^0||_{G^{-1}}^2,
\end{equation}
where $\rho:= || I - G \mathbb{E}[V] ||_{G^{-1}} = 1 - \lambda_{\min} (G^\frac12 \mathbb{E}[V] G^\frac12 ) <1$.

The result \eqref{eqn:K6convergence} is a general convergence for the type-\textsc{k} methods, since K1-5 can be regarded as the special cases of \ref{eqn:tagK6}. In the follwing we will investigate the convergence of the type-\textsc{c} iterative schemes following the techniques in \cite{GowerRichtarik_SIMAX15}.
Define $W:= \Omega (\Omega^T A^T G A \Omega)^{\dag} \Omega^T $, and $\widehat{G} := A^T G A $.  Note that $G$ is \textsc{spd}, then $\widehat{G}$ is also \textsc{spd} under the assumption that $A$ is of full column rank. It is easy to verify that $\Xi A = W \widehat{G}$, and the error propagation of scheme \ref{eqn:tagC6}  reads $e^{k+1} = (I-W \widehat{G}) e^k$. Taking expectation two times similar to \cite[Theorem 4.4]{GowerRichtarik_SIMAX15}, we get
$$ \mathbb{E} [ e^{k+1} ] = ( I - \mathbb{E} [W] \widehat{G} )\mathbb{E}[ e^k ].$$
Taking the $\widehat{G}$-norms to both sides, we have the estimate on the norm of expection
\begin{equation}
 || \mathbb{E} [ e^k ] ||_{\widehat{G}} \leqslant  \rho || \mathbb{E}[ e^{k-1} ] ||_{\widehat{G}} \leqslant \cdots \leqslant \rho^k || e^0 ||_{\widehat{G}},
\end{equation}
where $\rho = || I - \mathbb{E} [W] \widehat{G}||_{\widehat{G}} = 1- \lambda_{\min}(\widehat{G}^\frac12 \mathbb{E} [W] \widehat{G}^\frac12) = 1- \lambda_{\min}( \mathbb{E} [W] \widehat{G} )$.

Suppose that $\mathbb{E} [W]$ is positive definite. Similar to \cite[Lemma 4.5]{GowerRichtarik_SIMAX15}, we can prove that $\mathbb{E} [W] \geqslant (1-\rho) \widehat{G}^{-1}$. That is, $\forall ~ y \in \mathbb{R}^n$,
$$ \langle \mathbb{E} [W] \widehat{G} y, \widehat{G} y \rangle \geqslant (1-\rho) ||y||_{\widehat{G}}^2. $$

It is obvious to check that
$$ ||e^{k+1}||_{\widehat{G}}^2 = ||(I- W \widehat{G})e^k||_{\widehat{G}}^2 = ||e^k||_{\widehat{G}}^2 - \langle \mathbb{E} [W] \widehat{G} e^k, \widehat{G} e^k \rangle . $$
Taking expectation conditioned on $e^k$ yields
$$ \mathbb{E} [ ||e^{k+1}||_{\widehat{G}}^2 | e^k ] = ||e^k||_{\widehat{G}}^2 - \langle \mathbb{E} [W] \widehat{G} e^k, \widehat{G} e^k \rangle \leqslant \rho ||e^k||_{\widehat{G}}^2 . $$
Taking expectation again gives the following estimate on the expectation of norm.

\begin{theorem}
With the notations above, it holds that
$$\mathbb{E} [ ||e^{k+1}||_{\widehat{G}}^2 ]   \leqslant \rho  ~ \mathbb{E} [ ||e^k||_{\widehat{G}}^2 ].$$
Unrolling the recurrence we have the estimate
\begin{equation} \label{eqn:C6convergence}
 \mathbb{E} [ ||e^k||_{\widehat{G}}^2 ]   \leqslant \rho^k  ~   ||e^0||_{\widehat{G}}^2,
\end{equation}
where $\rho <1$ if $\mathbb{E} [W]$ is positive definite.
\end{theorem}
This result \eqref{eqn:C6convergence} has the similar form as the estimate \eqref{eqn:K6convergence} (see \cite[Theorem 4.6]{GowerRichtarik_SIMAX15}).

We consider the reasonable assumptions such that $\mathbb{E} [W]$ is positive definite.
Suppose that
\begin{enumerate}[(i)]
\item A random matrix $\Omega=\Omega_j \in \mathbb{R}^{n \times l_j}$ has a discrete distribution with probability $p_j >0$, such that $A \Omega_j$ ($j=1,\cdots, s$) has full column rank.
\item ${\bf \Omega} := [\Omega_1, \cdots, \Omega_s] \in \mathbb{R}^{n \times \sum_{j=1}^s l_j}$ has full row rank.
\end{enumerate}
Define the block diagonal matrix
$$ D = \text{diag} \left( p_1 (\Omega_1^T A^T G A \Omega_1)^{-1}, \cdots, p_s (\Omega_s^T A^T G A \Omega_s)^{-1} \right). $$
We can check that
\begin{equation}
\mathbb{E} [W] = \sum_{j=1}^s p_j \Omega_j (\Omega_j^T A^T G A \Omega_j)^{-1} \Omega_j^T = {\bf \Omega} D^\frac12 D^\frac12  {\bf \Omega}^T >0 .
\end{equation}

Similarly, we can prove the convergence of type-\textsc{s} methods, but we omit the details here.

\begin{corollary}
The convergence result of the scheme \ref{eqn:tagS4} reads
\begin{equation}
 \mathbb{E} [ ||e^k||_A^2 ]   \leqslant \rho^k  ~   ||e^0||_A^2,
\end{equation}
where $\rho = || I - \mathbb{E} [W] A||_A = 1- \lambda_{\min}(A^\frac12 \mathbb{E} [W] A^\frac12) $ with $W= \Omega (\Omega^T A \Omega)^{-1} \Omega^T $.
\end{corollary}

\subsection{Convergence analysis using random projections}

In the following we prove the convergence of type-\textsc{k} methods, and then extend the analysis to type-\textsc{c} and type-\textsc{s} methods.
%
Unless explicitly stated, in this subsection we consider the full-rank overdetermined problem, i.e., $m \geqslant n$ and $\mathrm{rank}(A) = n$.
We now consider the error propagation matrix $T$ of \ref{eqn:tagK6}  (see  Table \ref{tab:KCDseriesMethods})
\begin{equation}\label{eqn:TmatrixK6}
T = I_n - G A^T\Omega(\Omega^TA G A^T\Omega)^{-1}\Omega^T A,
\end{equation}
where $G$ is \textsc{spd}, and $\Omega=[\omega_{ij}]$  is an $m\times l$ real random Gaussian matrix with i.i.d. standard Gaussian normal random variables
$\omega_{ij}\sim N(0,1)$. We partition
$\Omega=[\omega_1,\ldots,\omega_l]$ into $l$ columns, where
$$
\omega_{j}=[\omega_{1j},\ldots,\omega_{mj}]^T, ~ j=1,\ldots, l.
$$
It is easily seen that the distribution density of $\Omega$ is given
by
\begin{eqnarray*}
p(\Omega) &=& \prod^{m,l}_{i=1,j=1}p(\omega_{ij}) =
\prod^{m,l}_{i=1,j=1}\frac1{\sqrt{2\pi}}\exp\big( -\frac12
\omega^2_{ij}\big)\\
&=&\frac1{\sqrt{(2\pi)^{ml}}}
\exp\big(-\frac12\mathrm{Tr}(\Omega^T\Omega)\big).
\end{eqnarray*}
 We have
$ \int p(\Omega)[d\Omega]=1$, where  $[d\Omega]=\prod^{m,l}_{i=1,j=1}d\omega_{ij}$ is the Lebesgue volume element.

Recall the vector-matrix correspondence $\mathrm{vec}$, which is
defined as follows: for $\Omega=\sum_{i,j}\omega_{ij}e_if^T_j$,
\begin{eqnarray*}
\mathrm{vec}(\Omega) = \sum^{m,l}_{i=1,j=1}\omega_{ij} f_j \otimes e_i,
\end{eqnarray*}
where $\{e_i:i=1,\ldots,m\}$ and $\{f_j:j=1,\ldots,l\}$ are standard
orthonormal bases for $\mathbb{R}^m$ and $\mathbb{R}^l$,
respectively. The notation $\otimes$ stands for Kronecker tensor
product. Clearly $\mathrm{vec}(\Omega)$ is a $ml$-dimensional real
vector. We also know that
$$
\langle\mathrm{vec}(X),\mathrm{vec}(Y)\rangle=(X,Y)_{HS},
$$
where  $(\cdot,\cdot)_{HS}$ is Hilbert-Schmidt inner product over the matrix space, defined by $(X,Y)_{HS}:=\mathrm{Tr}(X^TY)$.
Denote the Frobineus norm by $\|X\|_F=\sqrt{(X,X)_{HS}}$.

We can verify that
$$
\int
\frac{\mathrm{vec}(\Omega)\mathrm{vec}(\Omega)^T}{\langle\mathrm{vec}(\Omega),\mathrm{vec}(\Omega)\rangle}p(\Omega)[d\Omega]
= \frac1{lm}I_l \otimes I_m.
$$

Using the fact that $ \text{tr}_1 \left( \mathrm{vec} (X) \mathrm{vec}(Y)^T \right) = X Y^T $, where $\text{tr}_1$ denotes taking the trace over the first factor, we have
$$
\int \frac{\Omega\Omega^T}{\|\Omega\|^2_F}p(\Omega)[d\Omega] =
\frac1mI_m.
$$
Note that the matrix 2-norm is no more than the F-norm,
thus
$$
\|\Omega\|^2_F \geqslant \|\Omega\|^2_2 =
\|\Omega^T\Omega\|_2.
$$
Then we have the estimate
\begin{eqnarray*}
\int
\frac{\Omega\Omega^T}{\|\Omega^T\Omega\|_2}p(\Omega)[d\Omega] >
\frac1mI_m.
\end{eqnarray*}

Since
$$
\Omega^TA G A^T\Omega\leqslant
\lambda_{\max}\big(\Omega^TA G A^T\Omega\big)I_l\leqslant
\|\Omega^T\Omega\|_2 \lambda_{\max}(A G A^T)I_l,
$$
it follows that
\begin{eqnarray*}
A^T\Omega\big(\Omega^TA G A^T\Omega\big)^{-1}\Omega^TA\geqslant
\frac1{\lambda_{\max}(A G A^T)}\frac{A^T\Omega\Omega^TA}{\|\Omega^T\Omega\|_2}.
\end{eqnarray*}

Define $\widehat{T}:= G^{-\frac12} T G^\frac12 = I_n - G^\frac12 A^T\Omega(\Omega^TA G A^T\Omega)^{-1}\Omega^T A G^\frac12$, which is an orthogonal projector.
We can verify that
$$
\widehat{T} \leqslant I_n -
\frac1{\lambda_{\max}(A G A^T)}\frac{G^\frac12 A^T\Omega\Omega^TA G^\frac12}{\|\Omega^T\Omega\|_2}.
$$
Based on this, we obtain that
\begin{eqnarray*}
\mathbb{E}[\widehat{T}] \leqslant I_n -
\frac1{\lambda_{\max}(A G A^T)} G^\frac12 A^T\left(\int\frac{\Omega\Omega^T}{\|\Omega^T\Omega\|_2}p(\Omega)[d\Omega]\right)A G^\frac12
\end{eqnarray*}
Therefore we  obtain
\begin{eqnarray} \label{eqn:ConvergenceE[T]}
\mathbb{E}[\widehat{T}] \leqslant I_n -
\frac{G^\frac12 A^T A G^\frac12}{m\lambda_{\max}(A G A^T)}<I_n.
\end{eqnarray}

Let $T_j$ is the projector of the form \eqref{eqn:TmatrixK6} used in the $j$th ($j=0,1,\cdots, k$) iteration step of scheme \ref{eqn:tagK6}. Define $\widehat{T}_j= G^{-\frac12} T_j G^\frac12$ ($j=0,1,\cdots, k$).
In the following, we estimate the expectation:
$ \mathbb{E}_{T_0,\ldots,T_k} [  || T_k\cdots
T_0 e^0||_{G^{-1}}^2 ]$.

Let $ \widehat{e}^0 = G^{-\frac12} e^0  $.
Direct calculation indicates that
\begin{eqnarray*}
&&\mathbb{E}_{T_0,\ldots,T_k} [  || T_k\cdots
T_0 e^0||_{G^{-1}}^2 ] = \mathbb{E}_{\widehat{T}_0,\ldots,\widehat{T}_k} [  || \widehat{T}_k\cdots \widehat{T}_0 \widehat{e}^0||^2 ] \\
&&=\mathbb{E}_{\widehat{T}_0,\ldots,\widehat{T}_k} [ {\widehat{e}^{0T}} \widehat{T}_0\cdots \widehat{T}_k\cdots \widehat{T}_0 {\widehat{e}^0} ] = \mathbb{E}_{\widehat{T}_0,\ldots,\widehat{T}_{k-1}} [ {\widehat{e}^{0T}} \widehat{T}_0\cdots \mathbb{E}[\widehat{T}_k] \cdots \widehat{T}_0 {\widehat{e}^0} ] \\
&&\leqslant \mathbb{E}_{\widehat{T}_0,\ldots,\widehat{T}_{k-1}} [ {\widehat{e}^{0T}} \widehat{T}_0\cdots \widehat{T}_{k-1} \left( I_n - \frac{G^\frac12 A^T A G^\frac12 }{m \lambda_{\max}(A G A^T)} \right) \widehat{T}_{k-1}\cdots \widehat{T}_0 {\widehat{e}^0} ] \\
&& = \mathbb{E}_{\widehat{T}_0,\ldots,\widehat{T}_{k-1}} [ {\widehat{e}^{0T}} \widehat{T}_0\cdots \widehat{T}_{k-1}\cdots \widehat{T}_0 {\widehat{e}^0} ] \\ && \quad - \frac{1}{m \lambda_{\max}(A G A^T)} \mathbb{E}_{\widehat{T}_0,\ldots,\widehat{T}_{k-1}} [ {\widehat{e}^{0T}} \widehat{T}_0\cdots \widehat{T}_{k-1} G^\frac12 A^T A G^\frac12 \widehat{T}_{k-1} \cdots \widehat{T}_0 {\widehat{e}^0} ]
\end{eqnarray*}
Using the fact that $G^\frac12 A^T A G^\frac12 \geqslant \lambda_{\min}( G^\frac12 A^T A G^\frac12 ) I_n$, we get
\begin{eqnarray*}
\mathbb{E}_{\widehat{T}_0,\ldots,\widehat{T}_k} [\widehat{e}^{0 T} \widehat{T}_0\cdots \widehat{T}_k\cdots \widehat{T}_1 \widehat{e}^0 ]  &\leqslant&
\left( 1- \frac{1}{m \kappa} \right) \mathbb{E}_{\widehat{T}_0,\ldots,\widehat{T}_{k-1}} [ \widehat{e}^{0 T} \widehat{T}_0\cdots \widehat{T}_{k-1}  \cdots \widehat{T}_0 \widehat{e}^0 ]  \\
&\leqslant& \cdots \leqslant  \left( 1- \frac{1}{m \kappa} \right)^{k+1} ||e^0||_{G^{-1}} ,
\end{eqnarray*}
where $\kappa = \frac{\lambda_{\max}(G^\frac12 A^T A G^\frac12)}{\lambda_{\min}(G^\frac12 A^T A G^\frac12)}$
is the 2-norm condition number of $G^\frac12 A^T A G^\frac12$. If $G^{-1} = A^T A$, we have $\kappa =1$.

\begin{theorem}
With the notations above, we have the following estimate for the scheme \ref{eqn:tagK6}.
\begin{eqnarray}
\mathbb{E}_{T_0,\ldots,T_k} [  || T_k\cdots
T_0 e^0||_{G^{-1}}^2 ] \leqslant
 \rho ~ \mathbb{E}_{T_0,\ldots,T_{k-1}} [  || T_{k-1} \cdots
T_0 e^0||_{G^{-1}}^2 ],
\end{eqnarray}
and
\begin{eqnarray}
\mathbb{E}_{T_0,\ldots,T_k} [  || T_k\cdots
T_0 e^0||_{G^{-1}}^2 ]  \leqslant \rho^{k+1}  ||e^0||_{G^{-1}}^2 ,
\end{eqnarray}
where $\rho =   1- \frac{1}{m \kappa} < 1$.
\end{theorem}
Since $e^0$ is fixed, it follows that $e^{0T} G^{-1} e^0$ is
a fixed nonnegative number. Letting $k \to \infty$ and
$m,\kappa\geqslant 1$ be fixed, and using the fact that $e^k = T_k \cdots T_1 T_0 e^0$, we then have the convergence property
\begin{equation}
\mathbb{E}[||e^k||_{G^{-1}}^2] = \mathbb{E}_{T_0,\ldots,T_k} [ ||T_k\cdots T_0 e^0 ||_{G^{-1}}^2 ]  \leqslant \rho^{k+1} || e^0 ||_{G^{-1}}^2  \to 0 ~ (k \to \infty).
\end{equation}

{\bf Remarks}. The scheme \ref{eqn:tagK2} is the special of  \ref{eqn:tagK6} with $G=I$, and $\omega$ in \ref{eqn:tagK2} is corresponding to $\Omega  \in \mathbb{R}^{m \times l}$ with $l=1$.
Since $\omega$ can be rewritten as
$\omega = ||\omega|| \cdot u ,$
where $u$ is a normalized vector and $||\omega||^2 $ has
the chi-squared distribution density, 
we have
$$
T =I_n - \frac{A^T \omega  \omega^T A}{\omega^T AA^T \omega} = I_n - \frac{A^T u  u^T A}{u^T AA^T u} ,
$$
and
$$
\mathbb{E}[T] = I_n - \int \frac{A^T u u^T A }{u^T A A^T u} d\mu(u),
$$
where $d\mu(u)$ is the uniform probability measure in the sense that $\int d\mu(u)=1$.
Using the facts that $\int  u u^T   d\mu(u) = \frac{1}{m} I_m$, where $d\mu(u) =\frac{\Gamma(\frac n2)}{\pi^{\frac n2}} \int \delta(1- \langle u,u \rangle)[du]$, and $\lambda_{\min}(AA^T) \leqslant u^T {AA^T} u \leqslant
\lambda_{\max}(AA^T)$, we again have \eqref{eqn:ConvergenceE[T]}.

We now consider the convergence of type-\textsc{c} methods.
The error propagation matrix $T$ of \ref{eqn:tagC6} (see  Table \ref{tab:KCDseriesMethods}) reads
\begin{equation}\label{eqn:TmatrixC6}
T = I_n - \Omega (\Omega^T A^T G A \Omega)^{-1}\Omega^T A^T G A =
 I_n - \Omega (\Omega^T \widehat{G} \Omega)^{-1}\Omega^T \widehat{G} ,
\end{equation}
where $\widehat{G}=A^T G A$, and $\Omega$ is an $n \times l$ real random Gaussian matrix.

Define $\widehat{T} = \widehat{G}^\frac12 T \widehat{G}^{-\frac12} =I_n - \widehat{G}^\frac12  \Omega (\Omega^T \widehat{G} \Omega)^{-1}\Omega^T \widehat{G}^\frac12  $, and we can estimate
\begin{equation}\label{eqn:ExpectationTinC6}
 \mathbb{E}[ \widehat{T} ] \leqslant I_n - {\widehat{G}}/ (n \lambda_{\max}(\widehat{G})) .
\end{equation}

Note that the iteration error $e^k = T_k \cdots T_0 e^0$, where $T_j$ ($j=0,1,\cdots,k$) is the error propagation matrix in scheme \ref{eqn:tagC6}. Defining $\widehat{T_j} = \widehat{G}^\frac12 T_j \widehat{G}^{-\frac12}$ and $\widehat{e}^0 = \widehat{G}^\frac12 e^0$, we can check that
$$ ||e^k||_{\widehat{G}}^2 = ||T_k \cdots T_0 e^0||_{\widehat{G}}^2 = ||\widehat{T}_k \cdots \widehat{T}_0 \widehat{e}^0||^2. $$

Taking expectation and using \eqref{eqn:ExpectationTinC6}, we have
\begin{eqnarray*}
\mathbb{E}_{T_0,\ldots,T_k} [  || T_k\cdots
T_0 e^0||_{\widehat{G}}^2 ] \leqslant
 \left(1-\frac{\lambda_{\min}(\widehat{G})}{n \lambda_{\max}(\widehat{G})}\right) ~ \mathbb{E}_{T_0,\ldots,T_{k-1}} [  || T_{k-1} \cdots
T_0 e^0||_{\widehat{G}}^2 ].
\end{eqnarray*}
We summarize the convergence estimate of type-\textsc{c} methods in the following theorem.

\begin{theorem}
The iteration error of scheme \ref{eqn:tagC6} has the estimate
\begin{eqnarray}
\mathbb{E} [  || e^k ||_{\widehat{G}}^2 ] \leqslant
 \rho^{k+1} ~   || e^0||_{\widehat{G}}^2 ],
\end{eqnarray}
where $\rho = 1-\frac{1}{n \kappa} < 1$ with $\kappa$ being the condition number of $A^T G A$.
\end{theorem}

For type-\textsc{s} methods, where the coefficient matrix $A$ is \textsc{spd} and of size $n \times n$,  we typically consider the scheme \ref{eqn:tagS4} with the error propagation matrix $T=I_n - \Omega (\Omega^T A \Omega)^{-1}\Omega^T A$, and we can prove the following convergence result. The derivation is similar, so we omit the details here.
\begin{corollary}
The iteration error of scheme \ref{eqn:tagS4} is estimated by
\begin{eqnarray}
\mathbb{E} [  || e^k ||_A^2 ] \leqslant  \rho^{k+1} ~   || e^0||_A^2 ],
\end{eqnarray}
where $\rho = 1-\frac{1}{n \kappa} < 1$ with $\kappa$ being the condition number of $A$.

\end{corollary}

\section{Numerical examples}

The numerical tests were intensively performed in \cite{GowerRichtarik_SIMAX15}. In this section we will give a few example to demonstrate the convergence behavior of the randomized iterative schemes in Table \ref{tab:KCDseriesSchemes}.

Suppose that the coefficient matrix $A$ is of size $m\times n$. We test the matrices \textsf{rand} and \textsf{sprandn} for type-\textsc{k} and type-\textsc{c} methods, and use \textsf{sprandsym} for testing type-\textsc{s} methods. 
Here \textsf{rand} is generated by the corresponding \textsc{Maltab} function, with the entries drawn from the standard uniform distribution on the open interval $(0,1)$. \textsf{sprandn} is generated by the \textsc{Matlab} sparse random matrix function sprandn($m$,$n$, density, rc), where density is the percentage of nonzeros and rc is the reciprocal of the condition number. In our tests, we set density$=1/ \log(mn)$, and rc$=1/\sqrt{mn}$ as \cite{GowerRichtarik_SIMAX15}.
%
\textsf{sprandsym} is a sparse \textsc{spd} matrix generated by sprandsym($n$, density, rc, type), where type=1.  
%
The exact solution $x^*$ is assumed to be a vector of all ones. Then the right hand side is determined by $b=A x^*$.

\begin{figure}[!htbp]\centering
\subfigure[][Residual vs. iteration ]
   {\includegraphics[scale=0.17]{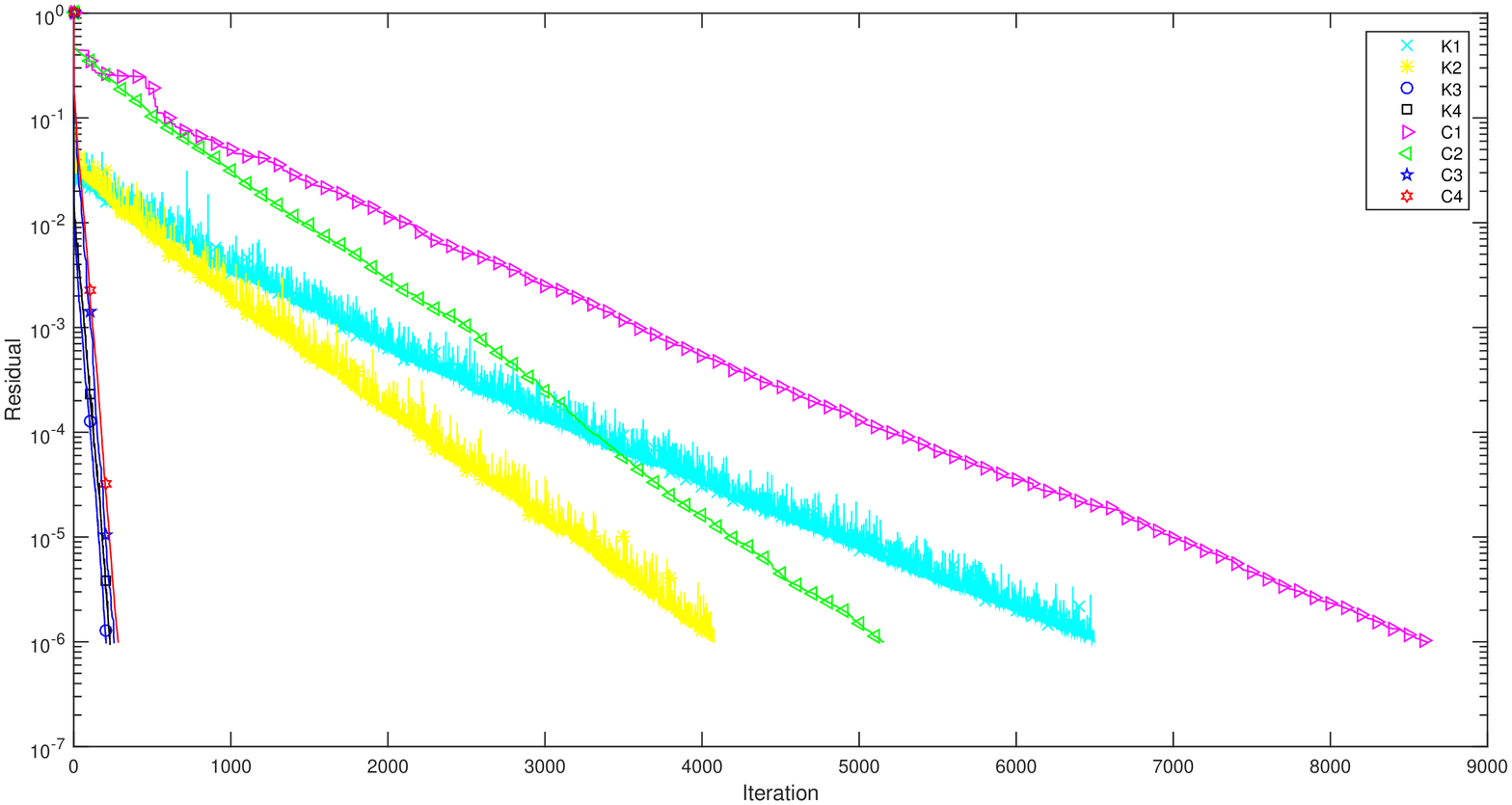}\label{fig:rand1}}
\subfigure[][Error vs. iteration ]
   {\includegraphics[scale=0.17]{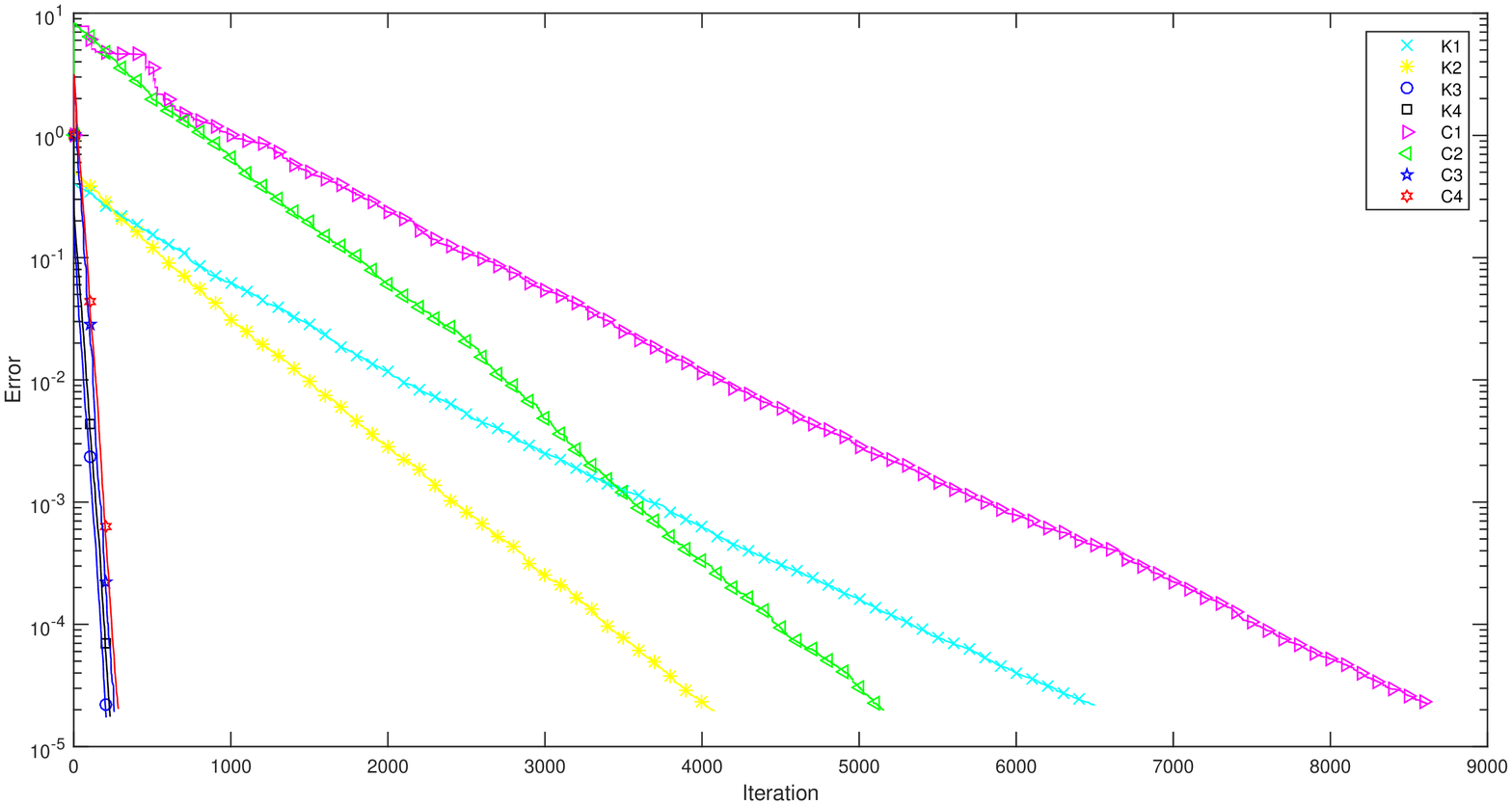}\label{fig:rand2}}
\subfigure[][Residual vs. time (s) ]
   {\includegraphics[scale=0.17]{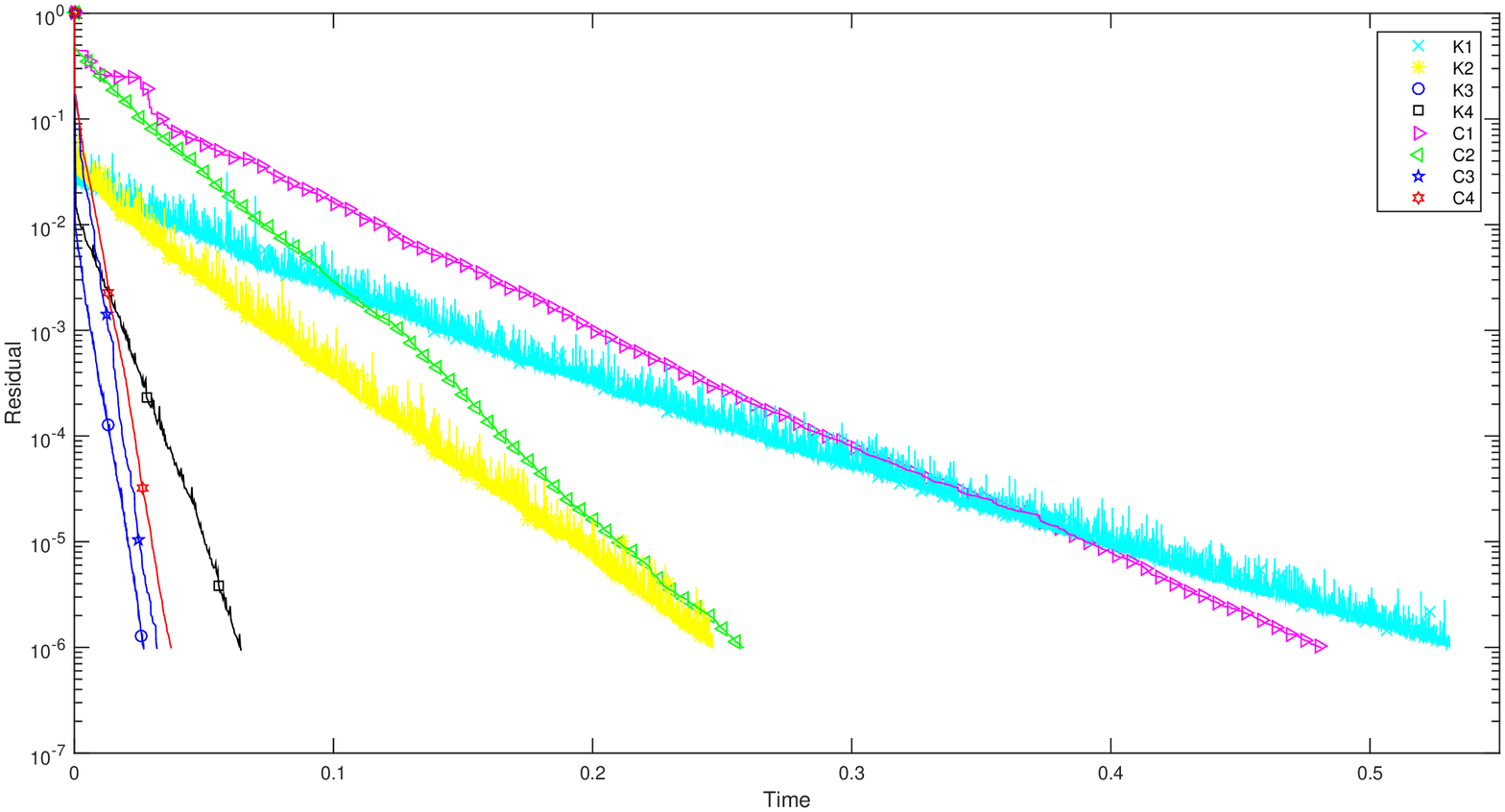}\label{fig:rand3}}
\subfigure[][Error vs. time (s) ]
   {\includegraphics[scale=0.17]{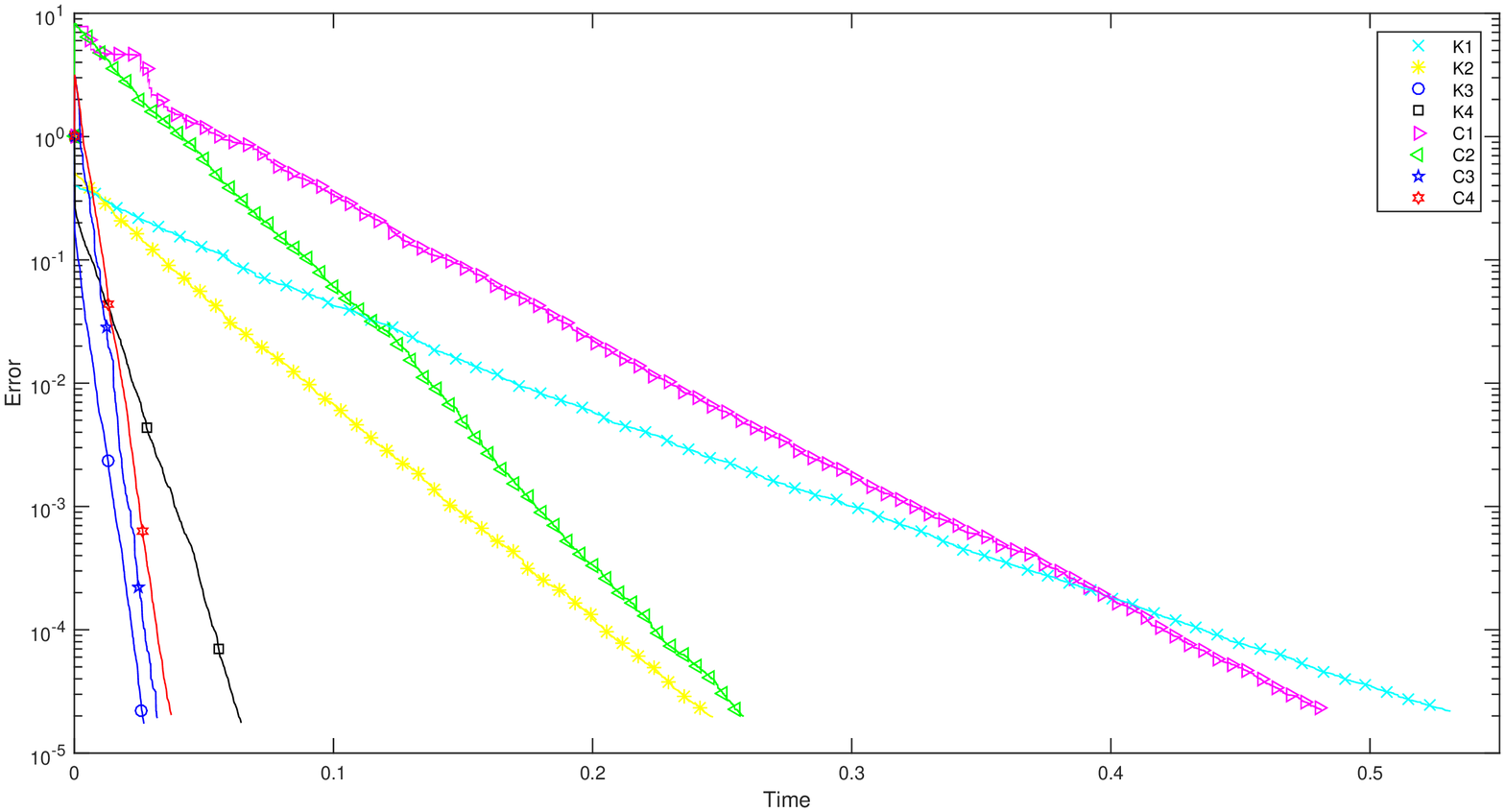}\label{fig:rand4}}
\caption{ \label{fig:rand} \textsc{rand} with $m = 1000$, $n = 100$.   }
\end{figure}

\begin{figure}[!htbp]\centering
\includegraphics[scale=0.35]{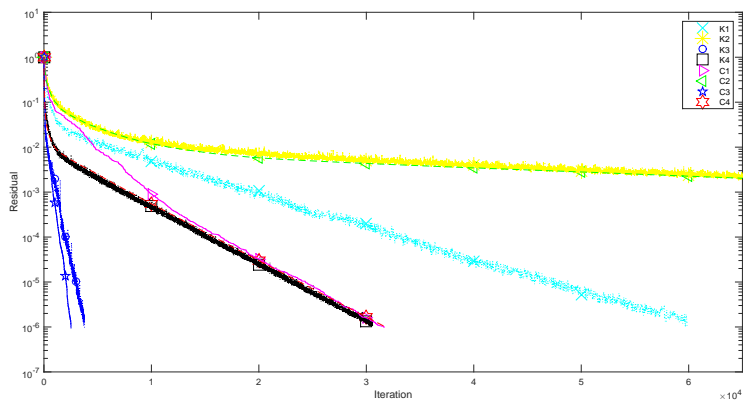} 
\caption{ \label{fig:sprandn} \textsc{sprandn} with $m = 100$, $n = 100$.    }
\end{figure}

\begin{figure}[!htbp]\centering
\includegraphics[scale=0.35]{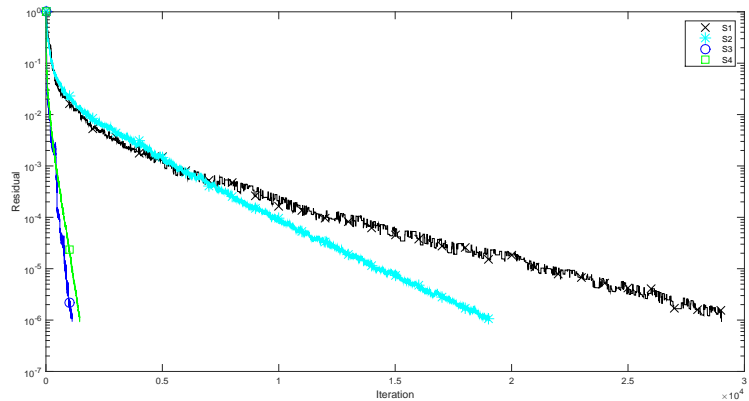} 
\caption{ \label{fig:sprandsym} \textsc{sprandsym} with $n = 100$.  }
\end{figure}

The the Schemes K5-6 and C5-6 depend on the choice of the matrix $G$. Such preconditioning techniques are another important and separate topic. For simplicity, we do not involve these here. We compare the schemes K1-4 and C1-4 for the unsymmetric cases (see Figures \ref{fig:rand} and \ref{fig:sprandn}), and D1-D4 for the \textsc{spd} case (see Figure \ref{fig:sprandsym}).
The performance of the iterative schemes also depend on the choice of the probability distribution \cite{GowerRichtarik_SIMAX15}. In this paper we do not focus on the probability distribution. 
For the schemes K1, K3, C1, C3, D1 and D3 using discrete sampling, we apply the \textsc{Matlab} function randsample($n$,$k$) to returns a $k$-by-1 vector of values sampled uniformly at random, without replacement, from the integers 1 to $n$.
For the schemes K2, K4, C2, C4, D2 and D4, we use randn($n$,$k$) to return an $n \times k$ Gaussian matrix with entries drawn from the standard normal distribution. For all the block version schemes, $k=$ floor(sqrt($n$)) as \cite{GowerRichtarik_SIMAX15}.

We choose the following parameters in the computation: the initial guess $x^0 = 0$, the maximum iteration number itmax = 100000, and the tolerance tol = 1e-6. When the iteration number is larger than itmax, or $||r||< \mathrm{tol} ||b||$, where the residual $r=b-Ax$, the iteration terminates. In each iteration step, we record the relative residual $\mathrm{res} = || b-Ax||/||b||$, and the relative error $\mathrm{err} = ||x-x^*|| / ||x^*||$. For each case, we record the wall-clock time measured using the tic-toc \textsc{Matlab} functions.

In Figures \ref{fig:rand1} and \ref{fig:rand2}, we
plot the relative residual and the relative error, respectively, on the vertical axis, and use the iteration number on the horizontal axis. The block version schemes, K3-4 and C3-4, need much less iteration steps. Even though in each step of the block version there exist matrix-matrix multiplication and a linear solver, the total computational time is much less; see Figures \ref{fig:rand3} and \ref{fig:rand4}, where the horizontal axis represents the computational time measured by using the tic-toc pair.
This is partly due to the fact that \textsc{Matlab} optimizes the matrix-matrix products and provides very efficient linear solvers.
This observation can be obtained in all test cases.

In Figure \ref{fig:rand},
for the single sample version iterative schemes as K1-2 and C1-2, the Gaussian methods 
require less iteration steps to reach a solution with the same precision as their discrete sampling counterparts. Despite the expensive matrix-vector product in each step required by the Gaussian methods, the computational time is also much less than the discrete sampling counterparts.
%
%
But in Figure \ref{fig:sprandn}, the Gaussian methods K2 and C2 need more iterations and more time.
%
%
We compare the four methods S1-4 on a system generated by the \textsc{Matlab} function \textsf{sprandsym}, for the single sample version scheme, S2 is faster than S1; the block version schemes S3-4 behave similarly, and are generally faster than the single sample ones. 

\subsection{Conclusion}

In this paper we present a unified framework to collect sixteen randomized iterative methods, such as Kaczmarz, \textsc{cd} and their variants. Under this general framework, we can recover the already known schemes, and derive three new iterative schemes as well. The convergence is proved under some general assumptions, for example, the coefficient matrix is of full column rank. But we believe that such restriction can be removed in the future work.
We give numerical examples to demonstrate the convergence behaviors of the iterative methods. The randomized strategies are as follows:
for methods based on discrete sampling we apply the uniform sampling without replacement, and for methods based on Gaussian sampling we use the Gaussian matrix with entries drawn from the
standard normal distribution. In this paper we do not focus on the probability distribution. But the choice of probability distribution can greatly affect the performance of the method and should be further investigated.

\end{document}